\documentclass[letterpaper, 10 pt, conference]{ieeeconf}  

\IEEEoverridecommandlockouts   

\usepackage{amsfonts}
\usepackage[all]{xy}
\usepackage{amsmath}
\usepackage{mathrsfs,amssymb,amsmath,booktabs,array,xcolor,url,color,soul,multirow}
\usepackage{bm}
\usepackage{arydshln}
\usepackage{enumitem}
\usepackage[marginal]{footmisc}

\usepackage{ntheorem}
\usepackage{graphicx}  
\usepackage{soul}
\soulregister\cite7
\soulregister\ref7 
\usepackage{wrapfig}
\usepackage{mathbbold}

\newtheorem{theorem}{Theorem}

\newtheorem{lemma}[theorem]{Lemma}

\newtheorem{example}{Example}
\newtheorem{definition}{Definition}

\title{\LARGE \bf
Capturing persistence of delayed complex balanced chemical reaction systems via decomposition of semilocking sets
}

\author{Xiaoyu Zhang and Chuanhou Gao* and Denis Dochain
\thanks{This work was funded by the National Nature Science Foundation of China under Grant No. 12071428 and 62111530247, and the Zhejiang Provincial Natural Science Foundation of China under Grant No. LZ20A010002.}
\thanks{X. Zhang is with Department of Control Science and Engineering, Zhejiang University, Hangzhou 310027, China
        {\tt\small Xiaoyu\_z@zju.edu.cn}}%
    \thanks{C. Gao is with School of Mathematical Sciences, Zhejiang University, Hangzhou 310027, China
        {\tt\small gaochou@zju.edu.cn (Corresponding author)}}%
\thanks{D. Dochain is with ICTEAM, UCLouvain, B\^{a}timent Euler, avenue Georges lema\^{i}tre 4-6, 1348 Louvain-la-Neuve, Belgium
        {\tt\small denis.dochain@uclouvain.be}}%
}    

\begin{document}
\maketitle
\thispagestyle{empty}
\pagestyle{empty}

\begin{abstract}
With the increasing complexity of time-delayed systems, the diversification of boundary types of chemical reaction systems poses a challenge for persistence analysis. This paper focuses on delayed complex balanced mass action systems (DeCBMAS) and derives that some boundaries of a DeCBMAS can not contain an $\omega$-limit point of some trajectory with positive initial point by using the method of semilocking set decomposition and the property of the facet, further expanding the range of persistence of delayed complex balanced systems. These findings demonstrate the effectiveness of semilocking set decomposition to address the complex boundaries and offer insights into the persistence analysis of delayed chemical reaction network systems.
\end{abstract}

\section{INTRODUCTION}
Chemical reaction networks (CRNs) are widely used in system modeling and analysis in various fields, including biochemical process \cite{Arceo2015}, electricity \cite{Samardzija1989}, medicine \cite{Allen2010} and even machine learning \cite{Anderson2021}. However real systems can be complex, which presents challenges for analysis. To address this, a common method of model reduction is to introduce time delays to replace complex intermediate processes \cite{G2019}. Time delays can also be used to induce gene switches in biological systems \cite{Wang2012}, model transport systems \cite{G2010}, and more. However, the introduction of time delays dramatically changes the dynamical properties of the system, making research on time-delayed chemical reaction networks essential.

The mass action law is a widely used approach for characterizing reaction rates and intensities in chemical reaction networks. Such a network with mass action kinetics is called a mass action system. Due to the complex interactions between multiple species in such systems, the dynamics of mass action systems can exhibit high levels of nonlinearity, posing challenges for their analysis. In response to this challenge, Chemical Reaction Network Theory (CRNT) has been developed over the past 40 years to investigate the relationship between the structure of a chemical reaction network and its dynamical properties, building on the pioneering work by Feinberg and Horn \cite{Feinberg1972Complex, Horn1972,Horn1974}. More recently, \cite{G2018} developed a dynamic model of a chemical reaction network with constant delays using the classical chain method described in \cite{Repin1965}, opening up new avenues of research on delayed CRNT \cite{Zhang2022}.

Persistence is one of the important topics of CRNT which refers to the non-extinction of species in a system and was first introduced in ecological systems. It plays a crucial role in the portrayal of biodiversity and system stability \cite{Tom2015}, characterizing the long-term dynamics of complex systems like animal populations \cite{Isbell2015}, the spread of infectious diseases \cite{Ha2004}, and biochemical reaction systems \cite{E2018}. For chemical reaction systems, persistence means that all species existing in the beginning will not be used up forever. From the perspective of dynamical systems and equations, persistence implies that the lower limit of the trajectory is greater than zero as time tends to infinity. Additionally, for bounded systems, we can use the relationship between the $\omega$-limit set and the boundary to inscribe it, making persistence an important theoretical research topic \cite{Zhang20222}.
Feinberg presented the ``persistence conjecture," stating that all weakly reversible mass action systems are persistent, making it another research hotspot. Another well-known conjecture is the ``Global Asymptotic Stability (GAC) Conjecture," which states that all complex balanced systems are persistent. Although there has been a great deal of work around these two hypotheses \cite{Anderson2008, Anderson2010, Angeli2007, G2014} for non-delayed chemical reaction systems, the persistence analysis of delayed systems is just beginning. \cite{H2019} derived that a conservative delayed chemical reaction system is persistent if there also exists a conservative relation between species in semilocking set. \cite{Zhang2021} foucses on the delayed complex balanced system and obtained that a DeCBMAS system is persistence if the boundary to which semilocking set corresponds is a vertex or a facet of stoichiometric compatibility class, thus deriving the persistence of 2d DeCBMAS directly.

In this paper, we focus on studying the persistence of DeCBMASs that have a stoichiometric subspace of dimension greater than 2d. In such systems, boundaries can be complex and extend beyond the facet or vertex cases. To address this complexity, we propose the semilocking set decomposition method as a potential tool.
Using this method, we investigate the boundaries $F_W$, where $W$ can be divided into several subsets $W^{(p)}$ such that each $F_{W^{(p)}}$ is a facet. We examine these boundaries from two perspectives: whether there are common elements between $W^{(p)}$ and the sufficient conditions of $F_W$ exhibiting no $\omega$-limit point of any trajectory starting from a positive point can be found.
Thus deriving the persistence of DeCBMASs with such semilocking boundaries $F_W$ directly. This method provides new possibilities for simplifying the persistence analysis of high-dimensional systems.
This paper is structured as follows: Section \ref{sec:2} provides preliminaries on chemical reaction mass action systems with time delays and the fundamental concepts of persistence. In Section \ref{sec:3}, we derive the persistence of some higher-dimensional delayed complex balanced systems. Finally, we conclude in Section \ref{sec:4}.
\hfill 
 
\hfill
~\\~\\
\noindent{\textbf{Mathematical Notation:}}\\
\rule[1ex]{\columnwidth}{0.8pt}
\begin{description}
\item[\hspace{-0.5em}{$\mathbb{R}^n;\mathbb{R}^n_{\geq 0};\mathbb{R}^n_{>0};\mathbb{Z}^n_{\geq 0}:$}] $n$-dimensional real space; $n$-dimensional non-negative real space; $n$-dimensional positive real space;$n$-dimensional non-negative integer space.
\item[\hspace{-0.5em}{$\bar{\mathscr{C}}_{+}, \mathscr{C}_{+}$}]: $\bar{\mathscr{C}}_{+}=C([-\tau,0];\mathbb{R}^{n}_{\geq 0}), \mathscr{C}_{+}=C([-\tau,0];\mathbb{R}^{n}_{>0})$ the non-negative, positive continuous function vectors defined on the interval $[-\tau,0]$, respectively.
\item[\hspace{-0.5em}{$x^{y_{\cdot i}}$}]: $x^{y_{\cdot i}}\triangleq\prod_{j=1}^{n}x_{j}^{y_{ji}}$, where $x,y_{\cdot i}\in\mathbb{R}^{n}$.
\item[\hspace{-0.5em}{$\mathbbold{0}_{n}$}]: An $n$-dimensional vector with each element to be zero.
\item[\hspace{-0.5em}{$0^{0}$}]: The result is defined by $1$.
\end{description}
\rule[1ex]{\columnwidth}{0.8pt}
\section{Preliminaries}\label{sec:2}
In this section, we shall provide some background information on delayed chemical reaction systems and persistence.

A chemical reaction network (CRN) $\mathcal{N}=(S,C,R)$ consists of a set of $n$ chemical species denoted by $S=\{X_1, X_2,...,X_n\}$ that participate in $r$ chemical reactions denoted by $R=\{R_1,..., R_r\}$. Each reaction $R_i$ is of the following form:
\begin{equation} \label{eq:1}
R_i: \quad\sum^{n}_{j=1}y_{ji}X_j\stackrel{}{\longrightarrow} \sum^{n}_{j=1}y'_{ji}X_j,
\end{equation}
where the stoichiometric coefficients $y_{ji}, y'_{ji}\in \mathbb{R}^n_{\geq 0}$ are non-negative integers,
 and the vectors $y_{.i}=(y_{1i},...,y_{ni})^\top$ and $y'_{.i}=(y'_{1i},...,y'_{ni})^\top$ are called reactant complex and the product complex in the complex set $C$.
 The sets of all the species, complexes, and reactions are denoted by $S$, $C$, and $R$, respectively, and together they uniquely determine the CRN $\mathcal{N}$.
The reaction vector $y'_{.i}-y_{.i}$ represents the change in concentrations of each species when the reaction $R_i$ takes place. All of the reaction vectors span a \textbf{stoichiometric subspace $\mathscr{S}$} of the network, defined as
\begin{equation}
\mathscr{S}={\rm span}\{y'_{.i}-y_{.i}|i=1,...,r\}.
\end{equation}
The orthogonal complement of the stoichiometric subspace $\mathscr{S}$ is denoted by $\mathscr{S}^\bot=\{a\in \mathbb{R}^n|a^\top y=0,{~\rm for~all~} y\in\mathscr{S}\}$.

The rate of each reaction $R_i$, satisfying the mass-action law, can be evaluated as
\begin{equation}
\delta_i(x)=k_ix^{y_{.i}}\triangleq k_i\prod^{n}_{j=1} x_{j}^{y_{ji}},
\end{equation}
where $x_j\in \mathbb{R}_{\geq0}$ represents the concentration of species $X_j(j=1,...,n)$, $x=(x_1,...,x_n)^\top$ represents the state, and the positive real number $k_i$ is the reaction rate constant.
The dynamics of a \textbf{mass-action system} can be expressed as
\begin{equation} \label{eq:2}
\dot {x}(t)=\sum^{r}_{i=1}k_ix(t)^{y_{.i}}(y'_{.i}-y_{.i}),\quad t\geq 0.
\end{equation}
We usually use a quadruple $(S,C,R,k)$ to capture a mass-action system.
The introduction of time delays has no effect on the structure of the network, but greatly affects the dynamical properties. Thus a delayed mass-action system shares the same stoichiometric subspace and equilibrium with the corresponding mass-action system, but has different dynamics and non-negative stoichiometric compatibility class from the latter.  \cite{G20182,G2018} made extensive studies on delayed mass-action systems.

The time delay in a chemical reaction can cause a lag in the generation of the product, while the consumption of the reactant occurs instantaneously. The dynamics of a \textbf{delayed mass-action system} can be described by the following equation \cite{G20182,G2018}:
\begin{equation} \label{eq:dde}
\dot{x}(t)=\sum^{r}_{i=1}k_i[x(t-\tau_i)^{y_i}y'_i-x(t)^{y_i}y_i],\quad t\geq0.
\end{equation}
where $\tau_i\geq0$ for $i=1,\cdots, r$ are constant time delays. A delayed mass-action system can be denoted as $\mathcal{M}=(S,C,R,k,\tau)$ where $S,C,R$ are the set of species, complex, and reactions respectively, and ${k},{\tau}$ are the vectors of reaction rate constants and time delays respectively.

The solution space of the delayed system \eqref{eq:dde} is $\bar{\mathscr{C}}_{+}$. When $\tau_i=0$ for $i=1,\cdots, r$, the system (\ref{eq:dde}) reduces to \eqref{eq:2}. Each trajectory of the system can only appear in a part of the solution space, and it cannot cover the entire solution space. Therefore \cite{G2018} proposed an equivalent class decomposition of the phase space $\bar{\mathscr{C}}$ called the \textbf{non-negative stoichiometric compatibility class}. Each stoichiometric compatibility class is a forward invariant set of trajectories, i.e. the trajectory $x^{\theta}$ starting from $\theta$ always stays in the stoichiometric compatibility class $\mathcal{P}_{\theta}$ containing $\theta$. The definition of $\mathcal{P}_{\theta}$ for the delayed system \eqref{eq:dde} is given by:
\begin{equation}\label{eq:scc}
\mathcal{P}_\theta=\{\psi\in \bar{\mathscr{C}}_+\;|\;c_a(\psi)=c_a(\theta)\; {\rm for\; all\;} a\in\mathscr{S}^\bot\},
\end{equation}
where the functional $c_a:\bar{\mathscr{C}}_+\rightarrow \mathbb{R}$ is defined by
 \begin{equation}\label{eq:zc0}
  \begin{split}
c_a(\psi)&=a^{\top}\biggl[\psi(0)+\sum^r_{i=1}\biggl(k_i\int^0_{-\tau_i}\psi(s)^{y_i}ds\biggr)y_i\biggr].\\
   \end{split} 
\end{equation}

A positive vector $\bar{x}\in \mathbb{R}^n_{>0}$ is called a \textbf{positive equilibrium} of $\mathcal{M}$ if it satisfies $\dot{\bar{x}}=0$.
A positive equilibrium $\bar{x}$ is called a \textbf{complex balanced equilibrium} if for any complex $\eta\in\mathbb{Z}^n_{\geq 0}$ in the network, it satisfies the following condition:
\begin{equation}
\sum_{i:~y_{.i}=\eta}k_i\bar{x}^{y_{.i}}=\sum_{i:~y'{.i}=\eta}k_i\bar{x}^{y{.i}},
\end{equation}
where $y_{.i},~y'_{.i},~k_i$ denote the reactant complex, the product complex, and the reaction rate constant of the $i$-th reaction, respectively.

A system is called a \textbf{complex balanced system} if each equilibrium is a complex balanced equilibrium. Complex balanced systems have been widely studied due to their good dynamical properties, such as the existence and uniqueness of a positive equilibrium in each positive stoichiometric compatibility class, and the local asymptotic stability of each positive equilibrium \cite{G20182}.

The following part presents the definition of persistence for DeMASs, which was first proposed by Komatsu and Nakajima \cite{H2019} and shares a similar form with that in \cite{Anderson2008} for the non-delayed case.

\begin{definition}[Persistence]
A DeMAS $(S,C,R,k,\tau)$ described by (\ref{eq:dde}) is persistent if any forward trajectory $x^{\psi}(t)\in\mathbb{R}^{n}{\geq 0}$ with a positive initial condition $\psi\in \mathscr{C}{+}$ satisfies
\begin{equation*}
\liminf_{t\rightarrow \infty}x^{\psi}_{j}(t)>0~~~\text{for all}~j\in{1,\cdots,n}.
\end{equation*}
\end{definition}
For bounded systems, persistence can be characterized by the $\omega$-limit set shown as follows.
\begin{definition}[$\omega$-limit set]
The $\omega$-limit set for the trajectory $x^{\psi}(t)$ with a positive initial condition $\psi\in \mathscr{C}_{+}$ is
\begin{equation*}
\begin{split}
    \omega(\psi):=&\{\phi\in \bar{\mathscr{C}}_{+}~|~x^{\psi}_{t_{N}}\rightarrow \phi,~\text{for some time sequence} \\&~t_{N}\rightarrow \infty~
    with ~t_{N}\in\mathbb{R}\}.
    \end{split}
\end{equation*}
The $\omega$-limit set is actually a positive invariant set of the corresponding trajectory.
\end{definition}

\begin{definition}[persistence for bounded trajectories]\label{def:2.3}
A DeMAS $(S,C,R,k,\tau)$ with the bounded trajectories is persistent if
\begin{equation}
 \omega(\psi)\cap(\cup_{W}L_W)=\varnothing
 , ~~~~\forall~\psi\in\mathscr{C}_{+},
 \end{equation} 
where 
\begin{equation}\label{def:Lw}
   L_W=\left\{w\in \bar{\mathscr{C}}_{+}\big|\substack{w_j(s)=0,~ X_j\in W,\\w_j(s)\neq 0,~X_j\notin W,}~~\forall s\in [-\tau_{\rm{max}},0]\right\},
\end{equation}
is called a boundary of $\bar{\mathscr{C}}_{+}$.Note that $\cup_{W}L_W$ is the set of all boundaries of $\bar{\mathscr{C}}_{+}$. Also, $F_W=L_W\cap D_{\psi}$ denotes the boundary of the stoichiometric compatibility class $D_{\psi}$.
\end{definition}
Among all the boundaries, a \textbf{facet} $F_W$ is a special case which satisfies $\dim{\mathscr{S}\vert_{W^c}}=\dim{\mathscr{S}}-1$.

The following concept plays an important role in the subsequent persistence analysis. 
\begin{definition}[semilocking set and locking set]
For a CRN $(S,C,R)$, a non-empty symbol set $W\subset S$ is called a semilocking set if it satisfies: $W\cap\mathrm{supp}~y_{\cdot i}\neq \varnothing$ if $W\cap\mathrm{supp}~y'_{\cdot i}\neq\varnothing$.
 If for any reaction $y_{\cdot i}\to y'_{\cdot i}$, there is $W\cap\mathrm{supp}~y_{\cdot i}\neq \varnothing$, $W$ is called a locking set.
\end{definition}
\section{MAIN RESULTS}\label{sec:3}
The goal of persistence analysis in chemical reaction networks is to determine whether the points on the semilocking boundary $F_W$ have the potential to be $\omega$-limit points of some trajectory with a positive initial point. Among all the boundaries, the ``facet" is a special case, and we have investigated whether $\omega$-limit points exist on facets for any DeCBMAS \cite{Zhang2021}. This result is formally stated in the following theorem.
\begin{theorem}[\cite{Zhang2021}]
Given a DeCBMAS $(S,C,R,k,\tau)$ of   (\ref{eq:dde}), for any semilocking set $W\subset S$, if $F_W$ defined by Definition \ref{def:2.3} is either empty or a facet of the stoichiometric subspace $\mathscr{S}$, this  DeCBMAS is persistent. 
\end{theorem}
As the complexity of chemical reaction networks increases, the dimension of the chemical stoichiometric space also increases, leading to a diversification of semilocking boundaries. In this paper, we conduct further analysis on the properties of other types of boundaries by decomposing the corresponding semilocking sets.

\subsection{Semilocking sets composed of independent subsets}\label{sec:3.1}
In this subsection, we focus on a semilocking set $W$ of a delayed chemical reaction mass action system $\mathcal{M}=\{S,C,R,k,\tau\}$ that can be decomposed into several independent subsets $W^{(p)}$, where each $F_{W^{(p)}}$ is a facet of $\mathcal{M}$. Here independent means that any two subsets have no common species, and any two species in different subsets have no interactions, i.e. they cannot participate in one reaction. We present an example to illustrate this case.
\begin{example}\label{ex:1}
Consider the following delayed system $\mathcal{M}$ 
\begin{equation*}
\begin{array}{c:c}
\begin{matrix}
 	\mathcal{M}:
  \xymatrix{2X_1\ar ^{\tau_1,~k_1~~} [r] &3X_{1}+X_{i} \ar ^{~\tau_2,~k_2} [d]\\
 & X_{1}+2X_i\ar ^{~\tau_3,~k_3~~} [lu]}\\
 \xymatrix{X_j+X_2\ar @{ -^{>}}^{~\tau_4,~k_4}@< 1pt> [r]& 2X_2 \ar  @{ -^{>}}^{~\tau_5,~k_5}  @< 1pt> [l]}\\
 \xymatrix{X_j\ar @{ -^{>}}^{~\tau_6,~k_6}@< 1pt> [r]& X_3 \ar  @{ -^{>}}^{~\tau_7,~k_7}  @< 1pt> [l]}
\end{matrix} &
\begin{matrix}
\mathcal{M}^{(1)}:~~~~~~&\\
\xymatrix{2X_1\ar ^{\tau_1,~k_1~~} [r] &3X_{1}+X_{i} \ar ^{~\tau_2,~k_2} [d]\\ 
& X_{1}+2X_i\ar ^{~\tau_3,~k_3~~} [lu]}\\
\mathcal{M}^{(2)}:~~~~~~&\\
 \xymatrix{X_j+X_2\ar @{ -^{>}}^{~\tau_4,~k_4}@< 1pt> [r]& 2X_2 \ar  @{ -^{>}}^{~\tau_5,~k_5}  @< 1pt> [l]}\\
     \xymatrix{X_j\ar @{ -^{>}}^{~\tau_6,~k_6}@< 1pt> [r]& X_3 \ar  @{ -^{>}}^{~\tau_7,~k_7}  @< 1pt> [l]}
\end{matrix}
\end{array}
\end{equation*}
$\mathcal{M}$ is a 4d weakly reversible network with zero deficiency, thus it is a DeCBMAS. And the subset $W=\{X_1,~X_2\}$ is a semilocking set with its boundary $L_W$ in the following form
\begin{equation*}
\begin{split}
     L_{W}=
 &\{\psi\in\bar{\mathscr{C}}_{+}\vert \psi (s)=(0,0,\psi_3(s),\psi_i(s),\psi_j(s))^\top, \\
 & s\in[-\tau_{\max},0],\psi_3(s)>0, \psi_i(s)>0,\psi_j(s)>0\}.
\end{split}
\end{equation*}
Vectors $(0,0,0,1,0)$ and $(0,0,1,0,-1)$ are both in the stoichiometric subspace $\mathscr{S}$. Thus $W$ is not a facet of $\mathcal{M}$. However, we can decompose it into two independent subsets $W^{(1)}={X_1}$ and $W^{(2)}={X_2}$ with $X_1, X_2$ not involved in the same reaction. Moreover, the boundaries $F_{W^{(1)}}$ and $F_{W^{(2)}}$ are both facets of $\mathcal{M}$.

There are two cases to consider:
\begin{enumerate}
\item $X_i$ and $X_j$ are not the same species. In this case, the subsets $\mathcal{M}^{(1)}$ and $\mathcal{M}^{(2)}$ are also independent, as there are no common species between them. Thus, the trajectory of each species in $\mathcal{M}$ is the same as that in $\mathcal{M}^{(1)}$ or $\mathcal{M}^{(2)}$. Hence the result that there does not exist an $\omega$-limit point of each trajectory with a positive initial point on $F_W$ is obvious.
\item $X_i$ and $X_j$ are the same species. In this case, there is coupling between the dynamics of $\mathcal{M}^{(1)}$ and $\mathcal{M}^{(2)}$. Therefore, we need a new approach to support further research.
\end{enumerate}
\end{example}
To address the coupling between the subsystems, one approach is to introduce a reduced system.
\begin{definition}[Reduced system]\label{def:re}
Let $\mathcal{M}=(S,C,R,k,\tau)$ be a delayed mass action system with dynamics $\dot{x}$, and let $\Tilde{S}\subset S$ be a subset of the species. The system $\Tilde{\mathcal{M}}=(\Tilde{S},\Tilde{C},\Tilde{R},\tilde{k}(t),\tilde{\tau})$ is called a reduced system of $\mathcal{M}$ if its dynamics can be expressed as $\dot{\Tilde{x}}=\dot{x}\vert_{\Tilde{S}}$, where $\dot{x}\vert_{\Tilde{S}}$ is the vector obtained by restricting the dynamics $\dot{x}$ to the species in $\Tilde{S}$. Thus each reaction of $\tilde{\mathcal{M}}$ has the form  $\Tilde{y}_{.i}\xrightarrow{\tilde{k}_i(t),~\tilde{\tau}_i} \Tilde{y}'_{.i}$ where 
$$\Tilde{y}_{.i}=y_{.i}\vert_{\Tilde{{S}}},~\Tilde{y}'_{.i}=y'_{.i}\vert _{\Tilde{S}}, ~\tilde{k}_i(t)=k_i\prod_{X_j\notin \tilde{S}}x_j^{y_{ji}},~\tilde{\tau}_i=\tau_i.$$
\end{definition}
According to above definition and setting $\Tilde{S}=W$, the delayed system $\mathcal{M}$ of the second case in Example \ref{ex:1} can be reduced as 
\begin{equation}
\begin{matrix}
 	\tilde{\mathcal{M}}:
  \xymatrix{2X_1\ar ^{\tau_1,~\tilde{k}_1(t)~~} [r] &3X_{1} \ar ^{~\tau_2,~\tilde{k}_2(t)} [d]\\
 & X_{1}\ar ^{~\tau_3,~\tilde{k}_3(t)~~} [lu]}\\
 \xymatrix{X_2\ar @{ -^{>}}^{\tau_4,~\tilde{k}_4(t)}@< 1pt> [r]& 2X_2(t) \ar  @{ -^{>}}^{\tau_5,~\tilde{k}_5(t)}  @< 1pt> [l]}
\end{matrix} 
\end{equation}
where $k_1(t)=k_1,~k_2(t)=k_2x_i(t),~k_3(t)=k_3x_i^2(t),~k_4(t)=k_4x_i(t),~k_5(t)=k_5$. Through reducing, not only $W$ but the system $\tilde{\mathcal{M}}$ can be divided into two independent parts. 
\begin{lemma}\label{lem:cf}
Let $\mathcal{M}=\{S,C,R,k,\tau\}$ be a DeCBMAS system and $W$ be a semilocking set of $\mathcal{M}$ that can be divided into independent subsets $W^{(p)}, p=1,\cdots, m$. Suppose that each $F_{W^{(p)}}$ is a facet of $\mathcal{M}$. Then an $\omega$-limit point cannot exist on the boundary of $F_W$.
\end{lemma}
\textbf{\textit{Proof.}}
As $W^{(p)}$ are independent, $\mathcal{M}$ can be divided into several subsystems $M^{(p)}$ according to $W^{(p)}$. And through reducing the system $\mathcal{M}$ based on the semilocking set $W$, each subsystem $\tilde{\mathcal{M}}^{(p)}$ of the reduced system $\tilde{\mathcal{M}}$ are also independent (no common species in $W^{(p)}$). Thus the trajectory of each species in $W$ of the original system $\mathcal{M}$ is the same as that of the reduced subsystem $\tilde{\mathcal{M}^{(p)}}$. So we just need to consider whether the origin of $\Tilde{M}^{(p)}$ can be an $\omega$-limit point of some trajectory. 
\\
$F_{W^{(p)}}$ is a facet of $\mathcal{M}$, then $\dim{\mathscr{S}\vert_{W^{(p)}}}=1$. Thus each reduced subsystem $\Tilde{\mathcal{M}^{(p)}}$ is a 1d complex balanced network with generalized mass action kinetics. And $v^{(p)}$ denotes the basis of the stoichiometric subspace of $\tilde{\mathcal{M}^{(p)}}$. \\
(1) If $v^{(p)}_j=0$ for some species $X_j$, the concentration of this species is a constant, namely, $x_j(t)=x_j(0)>0$ forever. So in this case, the origin can not be an $\omega$-limit point obviously.\\
(2) If $v^{(p)}_{j_1}\cdot v^{(p)}_{j_2}<0$ for some species $X_{j_1}, X_{j_2}$, for each positive initial point $\psi(s)$ of $\tilde{\mathcal{M}}^{(p)}$, each element in vector
$$\psi(0)+\sum^r_{i=1}\biggl(\int^0_{-\tau_i}k_i(s)\psi(s)^{y_i}ds\biggr)y_i-\mathbbold{0}_{\tilde{n}^{(p)}}$$ is positive, then it can not expressed by $v^{(p)}$. Thus in this case, the origin is not a point in the stoichiometric compatibility class of any positive initial point. Then the corresponding $F_{W^{(p)}}$ of the original system $\mathcal{M}$ is empty.
\\
(3) If each $v^{(p)}_j>0$, for each linkage class $L_l$ of $\tilde{\mathcal{M}}$, we can find one complex $\bar{y}^{l}\in L_l$ such that $\bar{y}^{l}_j-\tilde{y}_j<0$ for all $X_j\in \tilde{\mathcal{S}}$ and all complex $\tilde{y}\in L_l$.
We assume that the origin is an $\omega$-limit point of some trajectory $x^{\psi}(t)$ of $\mathcal{M}^{(p)}$ with positive initial point $\psi(s), s\in [-\tau,~0]$. 
Then for each $\epsilon>0$ and $t_0>0$, there exists some $t$ such that $x(t)$ is in the $\epsilon$-neighbourhood of the origin. Further combining the fact that 
\begin{equation}\label{eq:inf}
\lim\limits_{\tilde{x}\to 0}\frac{x^{\tilde{y}^{min}}}{x^{\tilde{y}}}=+\infty,
\end{equation}
we can obtain that for each constant $k>0$, there exists a small enough $\epsilon>0$, such that $x^{\Tilde{y}^{min}}>kx^{y}$ for each $x$ in the $\epsilon$-neighbourhood of zero.
Thus once the trajectory $x^\psi(t)$ enters into the $\epsilon$-neighbourhood of zero,  the dynamics each species $X_j\in \mathcal{S}$ can be expressed
\begin{equation}
\begin{split}
    \dot{x}_j&=\sum_{i=1}^r k_i(t-\tau_i)x^{\tilde{y}}(t-\tau_i)\tilde{y}_j'-\sum_{i=1}^r k_i(t)x^{\tilde{y}}(t)\tilde{y}_j\\
    &=\sum_{L_l}\sum_{i=1}^{\tilde{r}_l}[k_i(t-\tau_i)x^{\tilde{y}}(t-\tau_i)\tilde{y}_j'-k_i(t)x^{\tilde{y}}(t)\tilde{y}_j]
\end{split}
\end{equation}
From (\ref{eq:inf}), there exists an $\epsilon_1$ such that the sign of $\dot{x}_j$ is determined by
\begin{equation*}
\begin{split}
&\sum_{L_l}\sum_{\Tilde{y}_{.i}=\bar{y}^{l}}[k_i(t-\tau_i)x^{\bar{y}^l}(t-\tau_i)\tilde{y}'_{ij}-k_i(t)x^{\bar{y}^l}(t)\bar{y}^l_j]\\
&=\sum_{L_l}\sum_{\tilde{y}_{.i}=\Bar{y}^l}k_i(t)x^{\Bar{y}^l}(t)\tilde{y}'_{ij}(\frac{k_i(t-\tau_i)x^{\bar{y}^l}(t-\tau_i)}{k_i(t)x^{\Bar{y}^l}(t)}-\frac{\Bar{y}_j^{l}}{\tilde{y}'_{ij}})
\end{split}
\end{equation*}
As $\tilde{y}'_{ij}>\bar{y}^l_j>0$, $\frac{\Bar{y}_j^{l}}{\tilde{y}'_{ij}}<1$. Further combining 
\begin{equation*}
\lim\limits_{\tilde{x}\to 0}\frac{k_i(t-\tau_i)x^{\bar{y}^l}(t-\tau_i)}{k_i(t)x^{\Bar{y}^l}(t)}=1,
\end{equation*}
there exists an $\epsilon_2>0$ such that 
$$
\frac{k_i(t-\tau_i)x^{\bar{y}^l}(t-\tau_i)}{k_i(t)x^{\Bar{y}^l}(t)}-\frac{\Bar{y}_j^{l}}{\tilde{y}'_{ij}}>0,~\mathrm{for~each~}L_l.$$
So once the trajectory $x^\psi(t)$ enters into the $\epsilon$-neighbourhood of zero where $\epsilon=\min\{\epsilon_1,\epsilon_2\}$, there exists $x_j(t)>0$ for some $X_j$. This is obvious contradict with the fact that the origin is an $\omega$-limit point of the trajectory $x^\psi(t)$. 

 We can conclude that $F_W$ can not contain an $\omega$-limit point of any trajectory with positive initial point. 
$\hfill\blacksquare$

Thus the results of the persistence of DeCBMAS can be generalized by using Lemma \ref{lem:cf}.
\begin{theorem}
Consider a DeCBMAS $\mathcal{M}=(S,C,R,k,\tau)$ of dynamics of the form (\ref{eq:dde})
and let $W\subset S$ be a semilocking set that can be divided into several independent subsets $W^{(p)}$, such that for each $p$, the set $F_{W^{(p)}}$ defined by Definition \ref{def:2.3} is either empty or a facet of the stoichiometric compatibility class.
Then the DeCBMAS $\mathcal{M}$ is persistent. 
\end{theorem}
\textbf{\textit{Proof.}}
The result can be concluded directly by using Lemma \ref{lem:cf}.
$\hfill\blacksquare$

Thus going back to Example 1, the system $\mathcal{M}$ is persistent regardless of whether $X_i$ and $X_j$ are the same species or not.
\subsection{Semilocking sets composed of subsets sharing with common species}
In this subsection, we consider the semilocking set $W$ which can be divided into several $W^{(p)}$. According to $W^{(p)}$, the system $\mathcal{M}$ can be seen as a combination of sub-systems $\mathcal{M}^{(p)}=\{S^{(p)},C^{(p)},R^{(p)},k^{(p)},\tau^{(p)}\}$, where $p=1,\cdots,m$. Each boundary $F_{W^{(p)}}$ is a facet of the stoichiometric compatibility class of the system $\mathcal{M}^{(p)}$. However, different from Subsection \ref{sec:3.1}, $W^{(p)}, p=1,\cdots, m$ can have common species here. Additionally, there should not exist any interaction between two species in $W$ that are not in the same $W^{(p)}$, namely, they cannot participate in the same reaction.

\begin{example}\label{ex:2}
    $\mathcal{M}$ is a delayed system in the following form which is slightly different with the system in Example \ref{ex:1}. 
  \begin{equation*}
\begin{array}{c:c}
\begin{matrix}
 	\mathcal{M}:~~~~~~&\\
  \xymatrix{2X_1\ar ^{\tau_1,~k_1~~} [r] &3X_{1}+X_{i} \ar ^{~\tau_2,~k_2} [d]\\
 & X_{1}+2X_i\ar ^{~\tau_3,~k_3~~} [lu]}\\
 \xymatrix{X_1+X_2\ar @{ -^{>}}^{~\tau_4,~k_4}@< 1pt> [r]& 2X_2 \ar  @{ -^{>}}^{~\tau_5,~k_5}  @< 1pt> [l]}\\ \xymatrix{X_i\ar @{ -^{>}}^{~\tau_6,~k_6}@< 1pt> [r]& X_3 \ar  @{ -^{>}}^{~\tau_7,~k_7}  @< 1pt> [l]}
\end{matrix} &
\begin{matrix}
\mathcal{M}^{(1)}:~~~~~~&\\
\xymatrix{2X_1\ar ^{\tau_1,~k_1~~} [r] &3X_{1}+X_{i} \ar ^{~\tau_2,~k_2} [d]\\ 
& X_{1}+2X_i\ar ^{~\tau_3,~k_3~~} [lu]}\\
\mathcal{M}^{(2)}:~~~~~~&\\
 \xymatrix{X_1+X_2\ar @{ -^{>}}^{~\tau_4,~k_4}@< 1pt> [r]& 2X_2 \ar  @{ -^{>}}^{~\tau_5,~k_5}  @< 1pt> [l]}\\
 \xymatrix{X_i\ar @{ -^{>}}^{~\tau_6,~k_6}@< 1pt> [r]& X_3 \ar  @{ -^{>}}^{~\tau_7,~k_7}  @< 1pt> [l]}
\end{matrix}
\end{array}
\end{equation*}
The system $\mathcal{M}$ is a 4-dimensional weakly reversible network with zero deficiency, making it a DeCBMAS. The subset $W=\{X_1, X_2\}$ is a semilocking set with a boundary $L_W$ that can be represented as follows:
\begin{equation*}
\begin{split}
     L_{W}=
 &\{\psi\in\bar{\mathscr{C}}_{+}\vert \psi (s)=(0,0,\psi_3(s),\psi_i(s))^\top, \\
 & s\in[-\tau_{\max},0],\psi_3(s)>0, \psi_i(s)>0\}
\end{split}
\end{equation*}
Both vectors $(0,0,1,0)$ and $(0,0,0,1)$ are in the stoichiometric subspace $\mathscr{S}$. Although $W$ is not a facet of $\mathcal{M}$, it can be decomposed into two subsets $W_1=\{X_1\}$ and $W_2=\{X_1, X_2\}$, where $F_{W_1}$ and $F_{W_2}$ are both facets of $\mathcal{M}^{(1)}$ and $\mathcal{M}^{(2)}$, respectively. Note that $W$ cannot be divided into the form in Example \ref{ex:1} because $X_1$ and $X_2$ participate in the same reaction. Therefore further analysis is needed to address this case.
\end{example}
By reducing the system $\mathcal{M}$ to $\Tilde{\mathcal{M}}$ based on the semilocking set $W$, as defined in Definition \ref{def:re}, we can address the coupling caused by the species $X_i$. However, the reduced subsystems $\mathcal{M}^{(1)}$ and $\mathcal{M}^{(2)}$ are not independent due to the common species $X_1$.

Fortunately, we can use the following result to handle this situation.
\begin{lemma}\label{lem:cf1}
Consider a DeCBMAS $\mathcal{M}=\{S,C,R,k,\tau\}$ and $W$ is a semilocking set of $\mathcal{M}$ that can be partitioned into subsets $W^{(p)},p=1,\cdots, m$. $\mathcal{M}$ can be divided into subsystems $\mathcal{M}^{(p)}=\{S^{(p)},C^{(p)},R^{(p)},k^{(p)},\tau^{(p)}\}$ accordingly. The corresponding boundary $F_W$ cannot exist an $\omega$-limit point of any trajectory with a positive initial point if the following conditions hold:
\begin{itemize}
\item Each $F_{W^{(p)}}$ is a facet of the subsystem $\mathcal{M}^{(p)}$.
\item $S^{(p)}\cap (W-W^{(p)})=\emptyset$.
\end{itemize}
\end{lemma}
\textbf{\textit{Proof.}} We start by reducing the system $\mathcal{M}$ based on $W$ according to Definition \ref{def:re} to derive the reduced system $\tilde{\mathcal{M}}$. The reduced system can be partitioned into $m$ 1d subsystems $\tilde{\mathcal{M}}^{(p)}=\{\Tilde{S}^{(p)},\tilde{C}^{(p)},\tilde{R}^{(p)},\Tilde{k}(t)^{(p)},\tilde{\tau}^{(p)}\}, p=1,\cdots,m$.

Assume $F_W$ exists an $\omega$-limit point of some trajectory with positive initial point of $\mathcal{M}$, then the origin will be the $\omega$-limit point of the reduced system $\tilde{\mathcal{M}}$. 
If $\tilde{S}^{(p_1)}\cap\tilde{S}^{(p_2)}=\emptyset$ for any $p_1$ and $p_2$, i.e., if $\mathcal{M}^{(p_1)}$ and $\mathcal{M}^{(p_2)}$ are independent, the situation reduces to Lemma \ref{lem:cf}.
If there exists a common species $X_j$ between the sets $W^{(p)}$, the reduced system $\tilde{\mathcal{M}}$ is composed of 1d generalized mass action systems that are not independent. The dynamics of species $X_j$ can be written as:
\begin{equation*}
  \dot{x}_j(t)=\sum_{p=1}^{m}\dot{x}^{(p)}_j(t)
  \end{equation*}
where $\dot{x}^{(p)}_j(t)$ denotes the dynamics of $X_j$ in generalized mass action subsystem $\mathcal{M}^{(p)}$ and $\dot{x}^{(p)}_j(t)=0$ if $X_j\notin W^{(p)}$. 
\\
(1) $X_j$ only exists in one $\tilde{\mathcal{M}}^{(p)}$. As $\tilde{\mathcal{M}}^{(p)}$ is 1d, similar with Lemma \ref{lem:cf}, we can find a $\bar{y}^{l}$ in each linkage class $L_l$ in $\mathcal{M}^{(p)}$ such that $\bar{y}^{l}_j-\Tilde{y}^{(p)}_j<0$ for each $X_j\in W^{(p)}$ and $\tilde{y}^{(p)}\in L_l$. Thus we can directly conclude that $\dot{x}_j(t)>0$ when $x(t)$ comes into the $\epsilon$-neighbourhood of zero for some $\epsilon>0$ from the proof of Lemma \ref{lem:cf}.\\
(2) $X_j$ is the common species which exists in more than one subsystems $\tilde{\mathcal{M}}^{(p)}$. Although there exists coupling between the dynamics $\dot{x}^{(p)}_j(t)$, it does not affect the value of the complex in each subsystem $\tilde{\mathcal{M}}$. Then as each subsystem is 1d, the minimal complex $\bar{y}^{l}$ for each linkage class and each $\mathcal{M}^{(p)}$ exists. Thus for each $p$, there exists a $\epsilon^p>0$ such that $\dot{x}_j^{(p)}(t)>0$ when $x(t)$ comes into the $\epsilon^p$-neighbourhood of zero. Then $\dot{x}_j(t)>0$ when the trajectory $x(t)$ comes into the $\epsilon$-neighbourhood of zero where $\epsilon=\min\{\epsilon^{p}\}$.

However, this contradicts the assumption that the origin is one $\omega$-limit point of some trajectory with a positive initial point of the reduced system $\mathcal{M}$. Thus we can conclude the result.
$\hfill\blacksquare$

Then by using above lemma, we can generalized the persistence of DeCBMASs.
\begin{theorem}
    Let $\mathcal{M}=\{S,C,R,k,\tau\}$ be a DeCBMAS consisting of several sub-DeCBMAS $\mathcal{M}^{(p)}=\{S^{(p)},C^{(p)},R^{(p)},k^{(p)},\tau^{(p)}\}$. If each semilocking set $W$ of $\mathcal{M}$ can be decomposed into several $W^{(p)}$ such that each $F_{W^{(p)}}$ is a facet of $\mathcal{M}^{(p)}$ with $S^{(p)}\cap (W-W^{(p)})=\emptyset$, then $\mathcal{M}$ is a persistent system.
\end{theorem}
\textbf{\textit{Proof.}} It is obviously from Lemma \ref{lem:cf1}.$\hfill\blacksquare$

\begin{figure}
    \center
    \includegraphics[height=6cm,width=8cm]{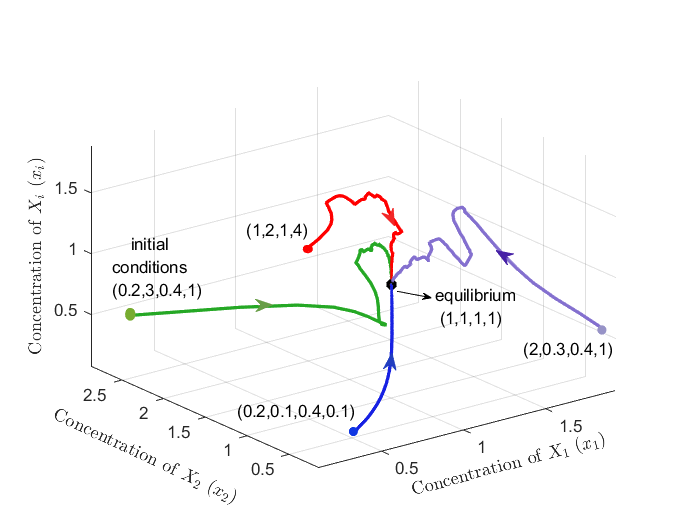}
    \caption{The evolution of the concentrations of $X_1,~X_2,~X_i$ with four different initial points of the system in Example \ref{ex:2}.}
    \label{fig:my_label}
\end{figure}

Based on the theorem stated above, we can conclude that the DeCBMAS $\mathcal{M}$ in Example \ref{ex:2} is a persistent system, as illustrated in Fig. \ref{fig:my_label}. The figure displays the evolution of the concentrations of the three species for four different initial points and delays. As we can see from the plot, all the trajectories converge to the unique equilibrium.


\section{CONCLUSIONS AND FUTURE WORKS}\label{sec:4}
The main focus of this paper is to generalize the persistence for DeCBMASs. Specifically, we aim to investigate whether there exist $\omega$-limit points of some trajectory starting from a positive initial point on a complex boundary that is not a facet or a vertex of the stoichiometric compatibility class. To achieve this goal, we focus on analyzing special complex boundaries $F_W$, where $W$ can be partitioned into several subsets $W^{(p)}$ such that each $W^{(p)}$ is a facet of a subsystem $\mathcal{M}^{(p)}$. By using this decomposition method and the properties of facets, we can determine whether the $\omega$-limit points are existence on these complex boundaries or not, thereby expanding the scope of DeCBMASs with persistence property.

In the future, we will give further consideration to the complex boundary $F_W$. We can divide $W$ into subsets $W^{(p)}$, and each subset $W^{(p)}$ will correspond to a facet or vertex of $F_{W}$. Additionally, we will consider cases where two species in the semilocking set $W$ but not in the same subset $W^{(p)}$ participate in the same reaction. This decomposition approach significantly simplifies the analysis of complex boundaries by breaking its corresponding semilocking set down into smaller, more manageable parts.






\begin{thebibliography}{99}
\bibitem{Arceo2015}
C. P. P. Arceo, E. C. Jose, A. Marin-Sanguino, and E. R. Mendoza,
Chemical reaction network approaches to Biochemical Systems Theory,
{\it Mathematical Biosciences}, vol. 269, 2015, pp 135-152.

\bibitem{Samardzija1989}
N. Samardzija, L.D. Greller, and E. Wasserman, Nonlinear chemical kinetic schemes derived from mechanical and electrical dynamical systems, {\it The Journal
of Chemical Physics}, vol. 90, no. 4, 1989, pp 2296–2304.

\bibitem{Allen2010}
L. Allen, An introduction to stochastic processes
with applications to biology, second edition. {\it Chapman
and Hall/CRC}, 2010.

\bibitem{Anderson2021}
D. F. Anderson, A. Deshpande, and B. Joshi, On
reaction network implementations of neural networks,
{\it Royal Society Interface}, vol. 18, pp 0–15.

\bibitem{G2019}
G. Lip\'{t}ak and K. M. Hangos, Distributed delay
model of the mckeithan’s network, {\it IFAC PapersOnline}, vol. 52, 2019, pp 33–38.

\bibitem{Wang2012} 
C. Wang, M. Yi, K. Yang, and L. Yang, Time
delay induced transition of gene switch and stochastic
resonance in a genetic transcriptional regulatory model.
{\it BMC Systems Biology}, vol.6, no. S9, 2012, pp 0-16.

\bibitem{G2010}
G. Orosz, R.E. Wilson, and G. Stépán, Traffic jams: dynamics and control, {\it Philosophical Transactions of the Royal
Society A}, vol. 368, 2010, pp 4455–4479.

\bibitem{Feinberg1972Complex}
M. Feinberg, Complex balancing in general kinetic  systems, {\it Archive for Rational Mechanics and Analysis},
vol.49, no. 3, 1972, pp 187–194.

\bibitem{Horn1972} 
F. Horn and R. Jackson, General mass action kinetics, {\it Archive for Rational Mechanics and Analysis}, vol. 47, no. 2, 1972, pp 81–116.

\bibitem{Horn1974}
M. Feinberg and F. Horn, Dynamics of open
chemical systems and the algebraic structure of the
underlying reaction network. {\it Chemical Engineering Science}, vol. 29, no. 3, 1974, pp 775–787. 

\bibitem{G2018}
G. Lip\'{t}ak, K.M. Hangos, and G. Szeder\'{k}enyi,
Approximation of delayed chemical reaction networks,
{\it Reaction Kinetics, Mechanisms and Catalysis}, vol. 123, no. 2, 2018, 403–419.

\bibitem{Repin1965}
Y. M. Repin, On the approximate replacement of
system with lag by ordinary dynamical systems, {\it Journal of Applied Mathematics and Mechanics}, vol. 29, pp 254–264. 

\bibitem{Zhang2022}
X. Zhang, C. Gao, and D. Dochain, On stability of two kinds of delayed chemical reaction networks, {\it IFAC PapersOnline}, vol. 55, no.18, 2022, pp 14-20.

\bibitem{Tom2015}
T. H. Oliver, M. S. Heard, N. J.B. Isaac et al., Biodiversity and resilience of ecosystem functions,
{\it Trends in Ecology \& Evolution}, vol. 30, Issue 11, 2015,
pp 673-684.

\bibitem{Isbell2015}
F. Isbell, D. Craven, J. Connolly et al., Biodiversity increases the resistance of ecosystem productivity to climate extremes, {\it Nature}, vol. 526, 2015, pp 574–577.

\bibitem{Ha2004}
T.J. Hagenaars, C.A. Donnelly, and N.M. Ferguson,
Spatial heterogeneity and the persistence of infectious diseases,
{\it Journal of Theoretical Biology}, Vol. 229, Issue 3, 2004,
pp 349-359.

\bibitem{E2018}
S. Ehrt, D. Schnappinger, and K. Y. Rhee, Metabolic principles of persistence and pathogenicity in Mycobacterium tuberculosis. {\it Nature Reviews Microbiology}, vol. 16, 2018, pp 496–507.

\bibitem{Zhang20222}
X. Zhang, Z. Fang, C. Gao, and D. Dochain, On the relation between $omega$-limit set and boundaries of mass-action chemical reaction networks, {\it Automatica}, vol 149, no. 110828, 2022. 

\bibitem{Anderson2008}
D. F. Anderson, Global Asymptotic Stability for A Class of Nonlinear Chemical Equations. {\it SIAM Journal on Applied Mathematics}, vol.68, no. 5, 2008, pp 1464–1476.

 \bibitem{Anderson2010}
D. F. Anderson and A. Shiu, The dynamics of weakly reversible population processes near facets, {\it SIAM Journal on Applied Mathematics}, vol. 70, no. 6, 2010, pp 1840–1858. 

\bibitem{Angeli2007}
D. Angeli, P. D. Leenheer, and E. D. Sontag, A petri net approach to the study of persistence in chemical reaction
networks, {\it Mathematical Biosciences}, vol. 210, no. 2, 2007, pp 98–618.

\bibitem{G2014}
M. Gopalkrishnan, E. Miller, and A. Shiu, A geometric approach to the global attractor conjecture. {\it SIAM Journal on Mathematical Analysis}, vol. 13, no. 2, 2014, pp 758–797. 

\bibitem{H2019}
H. Komatsu, and H. Nakajima, Persistence in
chemical reaction networks with arbitrary time delays.
{\it SIAM Journal on Applied Mathematics}, vol. 79, no. 1, 2019, pp 305–320.

\bibitem{Zhang2021}
X. Zhang and C. Gao, Persistence of delayed
complex balanced chemical reaction networks, {\it IEEE
Transactions on Automatic Control}, vol. 66, num. 4, 2021, pp 1658–1669.

\bibitem{G20182}
G. Lip\'{t}ak, K. M. Hangos, M. Pituk, and G. Szeder\'{k}enyi, Semistability of complex balanced kinetic
systems with arbitrary time delays. {\it System \& Control
Letters}, vol. 114, 2018, pp 38–43





\end{thebibliography}
\end{document}